\documentclass[12pt]{amsart}
\usepackage{amscd,amssymb}

\numberwithin{equation}{section}

\newtheorem {Theorem} 			{Theorem}
\newtheorem {varTheorem}                {Theorem}
\newenvironment {Theorem'}
        {\begin{varTheorem}{\hspace{-3.5mm}}{\bf '}{\hspace{3.5mm}}}
        {\end{varTheorem}}
\newtheorem {RefTheorem}[equation]     	{Theorem}  	
\newtheorem {Lemma}[equation]     	{Lemma}  	
\newtheorem {Proposition}[equation]	{Proposition}  

\newtheorem {Corollary}			{Corollary}

\theoremstyle{definition}

\newtheorem {Remark}[equation]		{Remark}

\newcommand{\pr} {\smallskip\noindent{\bf Proof\,\,}}

\newcommand     {\comment}[1]   {}
\newcommand     {\mute}[2] {}
\newcommand     {\printname}[1] {}

\newcommand{\labell}[1] {\label{#1}\printname{#1}}

\def	\red	{{\operatorname{red}}}

\def	\inv	{^{-1}}
\def	\to	{\longrightarrow}

\def	\R	{{\Bbb R}}

\def	\Q	{{\Bbb Q}}
\def	\Z	{{\Bbb Z}}
\def	\CP	{{\Bbb C}{\Bbb P}}

\def	\Mred	{{M_\red}}

\def	\pr	{\operatorname{pr}}

\begin{document}


\title[Semifree symplectic circle actions]
{On semifree symplectic circle actions
with isolated fixed points}

\author{Susan Tolman}
\author{Jonathan Weitsman}
\thanks{S. Tolman was partially supported by an NSF Mathematical Sciences
Postdoctoral Research Fellowship.
J. Weitsman was partially supported by NSF grant DMS 94/03567, by NSF 
Young Investigator grant DMS 94/57821, and
by an Alfred P. Sloan Foundation Fellowship.}

\address{Department of Mathematics,  University of Illinois at
Urbana-Champaign, 
Urbana, IL 61801}
\email{stolman@math.uiuc.edu}

\address{Department of Mathematics, University of California, Santa Cruz,
CA 95064}
\email{weitsman@cats.ucsc.edu}
\thanks{\today}

\begin{abstract}

Let $M$ be a symplectic manifold, equipped with a semifree symplectic circle
action with a finite, nonempty fixed point set.  We show that the circle
action must be Hamiltonian,
and $M$ must have the equivariant cohomology and Chern classes of $(P^1)^n$.

\end{abstract}

\maketitle



\section{Introduction}

Let $(M^{2n}, \omega)$ be a compact, connected symplectic manifold
of dimension $2n$.
A circle action on $M$  is {\bf symplectic}
if it preserves the symplectic form; that is, 
if its generating vector field $X$ satisfies $L_X \omega= d i_X\omega=  0$. 
A particular case is that of a Hamiltonian circle action,
where $i_X \omega$ is exact; in this case $i_X \omega = d \mu$ where
$\mu\in C^\infty (M)$ is the moment map.  A great deal is known
about Hamiltonian actions.  For example,  the quantization and
the push-forward measure  are determined by fixed point
data, and other manifold  invariants such as cohomology
and Chern classes are constrained  by this information.

A Hamiltonian circle action on a compact manifold
must have fixed points (one way to see this is that
the minimum of the moment map must be a fixed point).
In the case of a Kahler manifold
\cite{frankel}, or of a four-dimensional symplectic manifold \cite{mcduff},
the existence of a
fixed point guarantees that a symplectic circle action must in fact
be Hamiltonian.  This is not true in higher dimensions:  McDuff \cite{mcduff}
has constructed a symplectic six-manifold with a symplectic circle
action which has fixed points, but is not Hamiltonian. 
In this example the fixed point sets are tori; and there are no known
examples of symplectic, non-Hamiltonian circle actions with isolated
fixed points.  A natural question is whether any such examples can
exist.

In this paper we focus on the case of semi-free circle actions, and
show that any semi-free, symplectic circle action with
isolated fixed points is Hamiltonian
if and only if it has a fixed point.
Recall that an action of a group $G$ on a manifold $M$ is {\bf semi-free}
if the action is free on $M \setminus M^G$.\footnote{An action of a group 
$G$ on a 
manifold $M$ is {\bf quasi-free} if the stabilizers of points are connected;
in the case $G=S^1$ this is the same as the action being semi-free.}
Our main result is the following.

\begin{Theorem}
\labell{main}
Let $(M^{2n},\omega)$ be a compact, connected symplectic manifold, equipped
with a semifree, symplectic circle action with isolated fixed points.
Then if $M^{S^1}$ is nonempty, the circle action must be Hamiltonian.
\end{Theorem}

We prove this in Section 3 by an argument using integration
in equivariant cohomology. 

Theorem 1 brings us to the realm of Hamiltonian circle actions.  Though
a great deal is known about these, the general classification problem
for Hamiltonian $S^1$ spaces remains open.
However, they have been classified in
dimensions 2 and 4. 
In particular, up to symplectomorphism,  $P^1$
is the only two dimensional
example, and $(P^1)^2$
is the only four dimensional symplectic manifold 
with a semi-free Hamiltonian circle action with isolated fixed points.
In higher dimensions, much less is known: in the semi-free case with
isolated fixed points, the only example we know of 
is $(P^1)^n $, equipped with the diagonal circle action. 
As it turns out, the
classical manifold invariants---that
is, cohomology and Chern classes---of any such space must concord
with those of $(P^1)^n$.\footnote{The condition that
the fixed points be isolated implies the manifold is simply connected.}
Specifically, we prove the following theorem, which is an
equivariant version of Theorem \ref{Hattori} below,
due to Hattori \cite{hattori}.

\begin{Theorem}
\labell{coho}
Let $(M,\omega)$ be a compact, connected symplectic manifold, equipped
with a semifree, Hamiltonian circle action with isolated fixed points.
Let  $i:M^{S^1}\to M$, $j : {(P^1)^n}^{S^1}\to (P^1)^n$ denote the
natural inclusions of the fixed point sets of $M$ and $(P^1)^n$,
respectively.

There exists a map
from $M^{S^1}$ to $((P^1)^n)^{S^1}$ 
which identifies the  images of
$i^*:H^*_{S^1}(M,\Z)\to H^*_{S^1}(M^{S^1},\Z)$
and 
$j^*:H^*_{S^1}((P^1)^n,\Z)\to H^*_{S^1}(((P^1)^n)^{S^1},\Z)$.
This map sends the images of the equivariant
Chern classes of $M$ to those of $(P^1)^n$.
\end{Theorem}

For any  compact symplectic manifold $M$ with a Hamiltonian
circle action,   a theorem of Kirwan
\cite{kirwan} states that the natural inclusion map $i:M^{S^1}\to M$
of the fixed set $M^{S^1}$ into $M$ induces an injection
$i^*:H^*_{S^1}(M,\Q)\to H^*_{S^1}(M^{S^1},\Q)=H^*_{S^1}(M^{S^1},\Q)\otimes
H^*(BS^1,\Q)$.  In the case of a circle action with
isolated
fixed points, this injectivity theorem holds for integral
cohomology as well (see \cite{TW}).
Thus Theorem \ref{coho}  has the following corollaries:

\begin{RefTheorem}
There is an isomorphism
between $H_{S^1}^*(M,\Z)$ and $H_{S^1}^*((P^1)^n,\Z)$
which takes the equivariant Chern classes of $M$ to the
equivariant Chern classes of $(P^1)^n$.
\end{RefTheorem}

Thus the equivariant cohomology ring is
given by $H^*_{S^1}(M,\Z) =\Z[a_1,\ldots,a_n,y]/(a_iy-a_i^2)$,
and
the equivariant Chern series $c_t(M)= \sum_i t^i c_i(M)$ is given by
$$c_t(M) = \prod_i \left(1+t(2a_i - y)\right).$$

\begin{RefTheorem}(Hattori \cite{hattori})
\labell{Hattori}
There is an isomorphism
between $H^*(M,\Z)$ and $H^*((P^1)^n,\Z)$
which takes the Chern classes of $M$ to the Chern classes of
$(P^1)^n$.
\end{RefTheorem}

Thus, the cohomology ring is
given by $H^*(M,\Z) =\Z[a_1,\ldots,a_n]/(a_i^2)$,
and
the Chern series is given by
$$c_t(M) = \prod_i \left(1+t2a_i\right).$$

\begin{Remark} In dimension $6$, $M$ must be diffeomorphic to
$(P^1)^3$. This follows by a theorem of Wall \cite{Wall}
since  $M$ is a simply connected symplectic manifold with
the cohomology ring and Chern classes of $(P^1)^3$.
\end{Remark}

The basic idea underlying the proof of Theorem 2 is the use of 
Morse theory and integration in equivariant cohomology
to obtain a good picture of the cohomology ring of the
manifold in terms of the fixed point set.  This argument is given
in Section 4.

Our final result shows that the reduced spaces of any semifree Hamiltonian
$S^1$-space with isolated fixed points
have the cohomology rings and Chern classes of the
reduced spaces of $(P^1)^n$.  This corollary of Theorem 2 and of the
results of \cite{TW} is stated in Section 5 (Corollary 1).

\section{Equivariant Cohomology}

In this section we review  equivariant cohomology.
The integration formula (Proposition
\ref{integration}) is the main technical tool in this paper, and
will occur repeatedly.

Let $ES^1$ denote a contractible space on which $S^1$ acts
freely.
Let $BS^1 = ES^1/S^1$, and note
that $H^*(BS^1,\Z)$ is the polynomial ring in  a
single generator $x \in H^2(BS^1,\Z)$.

If $S^1$ acts on a manifold $M$, 
define $H^*_{S^1}(M,\Z) = H^*(M \times_{S^1} ES^1,\Z)$.  
In particular, if $P$ is a point, $H^*_{S^1}(P,\Z)$ is
naturally isomorphic to $H^*(BS^1,\Z)$; we will slightly
abuse notation by identifying these.

The projection $p: M \times_{S^1} ES^1 \to BS^1$ 
induces a pull-back map $H^*(BS^1,\Z)\to H^*_{S^1}(M,\Z)$;
this makes $H^*_{S^1}(M,\Z)$ into a $H^*(BS^1,\Z)$  module.
The projection $p$ also induces a
push forward map
$p_*:H^*_{S^1}(M,\Z)\to H^*(BS^1,\Z)$ given by 
``integration over the fibre''; we will usually denote $p_*$ by
the symbol $\int_M$.  
When $M$ is compact and the fixed points are isolated, this
has  the following expression:

\begin{Proposition}\labell{integration} Let 
$M$ be a compact manifold equipped with an action of
$S^1$ with isolated fixed points. 
Let $\alpha \in H^*_{S^1}(M,\Z)$.  Then as elements of $\Q(x)$,
$$\int_M \alpha  = \sum_F {\frac{\alpha|_{F} }{e(\nu_{F})}},$$

\noindent where the sum is taken over all
fixed points $F$, $\nu(F)$ is the normal bundle to $F$ and 
$e(\nu_{F})$ is the equivariant Euler class of this bundle.
\end{Proposition}

The right hand side of this equation is particularly simple
in the case where $M$ has an invariant almost complex
structure.
The computation of the
Chern classes of $\nu_F$ is then given by the following lemma.

\begin{Lemma}\labell{weights}
Let $E$ be a a representation of $S^1$, 
considered as complex $S^1$-equivariant vector bundle over a point $P$.
The representation $E$ decomposes as a direct sum
$E=\bigoplus E_i(w_i)$, where $E_i(w_i)$ is a complex line
on which the circle acts with weight $w_i$.
The equivariant Chern series $c_t(E)= \sum_i t^i c_i(E)$ is given by
$$c_t(E) = \prod_i(1+tw_ix),$$
\noindent where $x$ is the generator of $H^2_{S^1}(P,\Z)$.
\end{Lemma}

In particular, the Euler class of a representation
is $x^n$ times the product of the
weights of the circle action,
and the first Chern class is $x$ times the sum of the weights.

\section{Existence of the Moment Map: Proof of Theorem \ref{main}}
 
In this section, we prove Theorem \ref{main}.
We do this by  applying the integration formula of Proposition
\ref{integration} to the characteristic classes of our manifold to prove 
that if the manifold has any fixed points, then it has a fixed point
with no negative weights.  The following argument of
McDuff \cite{mcduff} then shows that this implies that the action is
Hamiltonian.

First, the symplectic form can be deformed to
a rational invariant symplectic form.
Furthermore, a symplectic manifold equipped with
a symplectic circle action
can be given a compatible invariant almost-complex structure.
Although this almost-complex structure is not unique,
both the Chern classes and the weights at every
fixed point are well-defined.
By Lemma \ref{count} below,
if $M$ has any fixed points
it has a fixed point $F$ with no negative weights.

In the case of a manifold with an invariant rational
symplectic form, the circle action will possess a circle valued 
moment-map.  Since  the fixed point $F$
is a local minimum for this circle-valued
moment map, it must lift to a real-valued moment map. 
This real-valued moment map is a Morse function with
only even index critical points, so
our manifold is simply connected. Therefore any symplectic
circle action on $M$ is Hamiltonian.

It remains to prove Lemma \ref{count}.

Note
that since the fixed points are all isolated, none of the weights
on the normal bundles to the fixed points can be zero.
Additionally, since the action is semifree, all these weights are $\pm 1$.

\begin{Lemma}\labell{count}
Let $M^{2n}$ be a compact manifold equipped with a semifree
circle action with isolated fixed points.  Suppose $M$ has an
invariant almost-complex structure.
Let $N_k$ denote the number of fixed points of the circle action
with $k$ negative weights. 
Then for all $k$ with $0 \leq k \leq n$
$$ N_k = N_0 \frac{n!}{ k! (n-k)!}.$$
\end{Lemma}

\begin{proof}
Let $P$ be a point, and let $x$ be the generator of
$H^2_{S^1}(P,\Z)$. 
Let $y \in H^2_{S^1}(M,\Z)$ be the  pull-back of $x$ to $M$ by the map
$p:M\to P$.
Consider the equivariant cohomology class $\gamma = \frac{1}{2}(ny - c_1(TM))$.

For dimensional reasons, for any $l$ with
$0 \leq l < n$,   
\begin{equation*}
\int_M \gamma^l = 0.
\end{equation*}
Applying the equivariant integration formula of Proposition
\ref{integration} yields
\begin{equation*}
\int_M \gamma^l = \sum_F {\frac { {\gamma^l}|_{F} }  {e(\nu_{F})} },
\end{equation*}
where $e(\nu_F)$ is the equivariant Euler class of the
tangent bundle at $F$ and the sum is over all fixed points $F$.

For every fixed point $F$, let $k_F$ be the number of negative weights
of the circle action on the normal bundle $\nu_F$.
Applying  Lemma \ref{weights}, we see that
$$\gamma|_{F} =  k_F \ x \ \ \ \mbox{    and    } \ \ \ 
e(\nu_{F}) = (-1)^{k_F} \  x^n.$$
Thus for all $l$ with $0 \leq l < n$,
\begin{equation*}
\sum_F (k_F)^l (-1)^{k_F} = 0.
\end{equation*}
Equivalently, for all $l$ with $0 \leq l < n$,
\begin{equation}
\labell{gamma}
\sum_{k=0}^n N_k k^l (-1)^k = 0.\end{equation}

Define a $n\times (n+1)$ dimensional matrix $V$ whose entries are given
by $$V_{i,j} = j^i$$
\noindent for $0 \leq i < n$ and $0 \leq j \leq n.$
Since this matrix consists of the first $n$ rows of
a nonsingular $(n+1) \times (n+1)$ Vandermonde matrix,
it has rank $n$.

Let  $A$ be the column vector whose entries are
$A_k=(-1)^k N_k$. 
Equation \ref{gamma} can be written
as a matrix equation $$V A =0.$$

Up to multiplication by a constant,
there is a unique solution to equation \ref{gamma}.  Moreover,  
the fixed point data of $(P^1)^n$ is a solution.
\end{proof}

\begin{Remark} Our methods generalize to the case of non-semifree circle
actions with
finite non-empty fixed point sets and
give constraints on the possible fixed point data of
such actions.
In the semifree case above these constraints eliminate all
non-Hamiltonian examples.  In the general
case, some possibilities cannot be ruled out by our methods.  For example,
we cannot rule out the existence of a symplectic six-manifold, equipped
with a symplectic circle action which has two fixed points, with the action
having weights $(1,1,-2)$ on the normal bundle to one fixed point
and $(-1,-1,2)$ on the normal bundle to the other.
\end{Remark}

\section{The cohomology ring: Proof of Theorem 2}

In this section we prove Theorem \ref{coho}.  The method
of proof again involves the equivariant fixed point formula
(Proposition \ref{integration}), but this time
we 
apply the formula not only to the Chern classes of the manifold
$M$ but also to additional
classes.
To construct these classes we  use the following theorems,
which are proved using  Morse theory.

\begin{RefTheorem}\labell{perfection}(Frankel, Kirwan)
Let a circle act on a 
a symplectic manifold $M$ with moment map $\mu: M \to \R$. 
Then the moment map $\mu$ is a perfect Morse function on $M$ (for both
ordinary and equivariant cohomology).
The critical points of $\mu$  are the fixed points of $M$,
and the index of a critical point $F$ is precisely twice the number 
of negative weights of the circle action on $TM_{F}$.
\end{RefTheorem}

\begin{RefTheorem}\labell{injectivity}(Kirwan)  
Let a circle act on a 
a symplectic manifold $M$ in a Hamiltonian fashion. 
Let $i :M^{S^1} \to M$ denote the natural inclusion.
The induced map $i^*:H^*_{S^1}(M^{S^1},\Q)\to H^*_{S^1}(M,\Q)$ is injective.
\end{RefTheorem}

More precisely, we need to take
the following result from Kirwan's proof of
Theorem \ref{perfection} and Theorem
\ref{injectivity} (\cite{kirwan}, see also \cite{TW}). 
For simplicity, we restrict our attention to the case
of circle actions with isolated fixed points.

\begin{RefTheorem}\labell{getclasses}
Let a circle act on a 
a symplectic manifold $M$ in a Hamiltonian fashion 
with isolated fixed points.
Let $F$ be any fixed point of index $2k$.
Let $w_1,\ldots,w_k$ be the negative weights of the circle action on $TM_{F}$.
Then there exists a class $a_F \in H^{2k}_{S^1}(M,\Z)$ such that
$a_F|_F = (-1)^k x^k \prod_{i=1}^k w_i$ and 
$a_F|_{F'} = 0$ for all fixed points $F'$ of index
less than  $2k$.
Moreover, taken together over all fixed points, these classes
are a basis for the cohomology  $H^*_{S^1}(M,\Z)$ as a $H^*(BS^1,\Z)$ module.

When the action is semi-free, there is a unique  way to choose 
the class
$a_F$ so that $a_F|_{F'} = 0$ for all other fixed points $F'$ of index
less than  or equal to $2k$.
\end{RefTheorem}

The first step is to analyze these cohomology classes on $(P^1)^n$.
Our strategy will then be to show their
counterparts on $M$ mimic their behavior. 
Each fixed point in $(P^1)^n$ is on the north pole
of some set of spheres, and on the south pole for all the other spheres.
Thus, if we identify the spheres with
the integers $1$ through $n$,
the fixed points $q_k^i$ of index $2k$ are in one-to-one
correspondence with subsets $J \subset \{1,\ldots,n\}$
with $k$ elements.  
Consider two fixed points $F$ and $F'$, 
which correspond to subsets  $J$ and $J'$ of $\{1,\ldots,n\}$,
with $k$ and $k'$ elements, respectively.
Let $\alpha_F$ be the cohomology class associated to $F$
as in Theorem \ref{getclasses}.
Then it is easy to see that
$\alpha_F|_{F'} =  x^k$ if $J \subset J'$, otherwise
$\alpha_F|_{F'} =  0$.
 
We now study the cohomology group $H^2_{S^1}(M,\Z)$. 
Denote by $p_k^i$ the critical points of index $2k$ in $M$.
By Theorem \ref{getclasses}, for
each $1 \leq j \leq n$, we can find a class $a_j$ such that
\begin{equation}\label{eq:4.1}
a_j|_{p_1^j} = x \quad \mbox{ and }\cr
a_j|_F = 0 \mbox { for all other critical points } F \mbox{ of index }
0 \mbox{ or } 1.
\end{equation}

In the remainder of this section, we will prove that these
forms satisfy the following Proposition: 

\begin{Proposition}\label{prop:4.5}
Let $J$ be a subset of $\{1,\dots,n\}$ with $k$ elements. 
There exists a unique fixed point $F$ of index $2k$ such
that $a_j|_F = x$ if and only if $j \in J$,
and  $a_j|_F = 0$ otherwise.
\end{Proposition}

By identifying the fixed points of
$M$ with subsets of $\{1,\ldots,n\}$,
this Proposition gives an isomorphism between
the fixed point set of $M$ and  the fixed point set
of $(P^1)^n$.
We claim that this map  
identifies the  images of
$i^*:H^*_{S^1}(M,\Z)\to H^*_{S^1}(M^{S^1},\Z)$
and $j^*:H^*_{S^1}((P^1)^n,\Z)\to H^*_{S^1}(((P^1)^n)^{S^1},\Z)$,
and also identifies the images of the equivariant Chern classes.

First, this map
sends the classes $a_i$ to the classes $\alpha_i := \alpha_{q_1^i}$. 
More generally, let $F$ be any fixed point of index  $2 k_F$, and let
$J_F$ be the corresponding subset of $\{1,\ldots,n\}$. 
The cohomology class $a_F = \prod_{j \in J_F} a_j$ is
the unique class which restricts to $x^k$  on $F$,
and to $0$ on all other fixed points of index less than
or equal to $2 k_F$.   Now for any
other fixed point $F'$ with an associated
subset $J_{F'}$ , $a_F|_{F'} = x^{k_F}$   exactly
if $J_{F} \subset J_{F'} $, and otherwise
$a_F|_{F'} = 0$.  Since the same  is true for the
$\alpha_F$, this map also  takes the classes $a_F$ to
the classes $\alpha_F$.  Since these classes form a 
basis for the cohomology as a $H^*(BS^1,\Z)$ module,
this proves the first claim. 

Moreover, the restriction of an equivariant
Chern class of any manifold to the fixed point set is
the equivariant Chern class of the normal bundle to
that fixed set.  For isolated fixed points,
this class is determined by the weights
of the circle action on the normal bundle.
Since this map takes every fixed point in $M$ to fixed point
in $(P^1)^n$ with the same weights, it takes
equivariant Chern classes to equivariant Chern classes.

Thus, it remains to prove
Proposition \ref{prop:4.5}. To do this we will repeatedly use the
following generalization of the
technique used in the last section.

\begin{Lemma}
\labell{technical}
Let $M$ be a compact symplectic manifold with a 
semi-free Hamiltonian circle action with isolated
fixed points.
Let $d \in H^{2l}_{S^1}(M,\Z)$  
and $\delta \in H_{S^1}^{2l}((P^1)^n,\Z)$
be  equivariant cohomology classes of degree $2l$.

If there exist  $l+1$ integers $k$ with
$0 \leq k \leq n$ such that
$$\sum_i d|_{p_k^i} = \sum_i \delta|_{q_k^i}, $$ 
then for {\bf all} integers $k$ with $0 \leq k \leq n,$
$$\sum_i d|_{p_k^i} = \sum_i \delta|_{q_k^i}. $$ 
Here, the sum is always taken over all fixed points of index $2k$.
\end{Lemma}

\begin{proof}
For dimensional reasons, for all $i$ with $0 \leq i <  n-l$,
\begin{equation*} 
\int_M d \cdot \gamma^i = 0 \end{equation*}
On the other hand, applying the formula in Proposition \ref{integration},
\begin{equation*} 
\int_M d \cdot \gamma^i = \sum_{k=0}^n \sum_i
{\frac{d|_{p_k^j} \ k^i \ x^i} {(-1)^k \  x^n}}. \end{equation*}
Thus, for all $i$ with $0 \leq i <  n-l$,
\begin{equation}\label{dform} 
\sum_{k = 0}^n 
\sum_i
 (-1)^k d|_{p_k^i} \, k^i = 0. \end{equation}

Define  a matrix $V^{n-l}$ which has
$n-l$ columns and $n+1$ rows, and with entries given
by $V^{n-l}_{ij} = j^i$,
where $0 \leq i \leq n$ and $0 \leq j < n-l$.
The matrix $V^{n-l}$ consists of  the first $n-l$ rows of an
$(n+1) \times (n+1)$ Vandermonde matrix.

Define a column vector $D$ by
$$D_k  = (-1)^k \sum_i d|_{p_k^i}$$ for
all $0 \leq k \leq n$,
where the sum is taken over all fixed points of index $2k$.
Then equation \ref{dform}
is equivalent to the matrix equation

\begin{equation*}  V^{n-l} D =0.\end{equation*}
Similarly, define a column vector $\Delta$ by
$\Delta_k  = (-1)^k \sum_i \delta|_{q_k^i}$ for
all $k$ with $0 \leq k \leq n$,
where the sum is taken over all fixed points of index $2k$.
An argument identical to that given above shows that 
the vector $\Delta$ satisfies the matrix equation
$ V^{n-l} \Delta =0.$

Finally, by assumption,
there are $l+1$ integers $0 \leq k \leq n$ such
that  $D_k = \Delta_k.$
The remaining entries of the difference $D-\Delta$
form an $(n-l)$-dimensional column vector which lies
in the kernel of a nonsingular $(n-l) \times (n-l)$-dimensional
Vandermonde matrix.  Thus $D=\Delta$.
\end{proof}

\begin{Lemma}\labell{aLemma} 
For any $1 \leq j \leq n$ and $0 \leq k \leq n$,
$$\sum_i a_j|_{p_k^i} = {\frac{(n-1)!}{(k-1)!(n-k)!}} \ x, $$
where the sum is taken over all fixed points of index $2k$.
In particular,  since there is only one fixed point of index $2n$,
\begin{equation}
\labell{atop}
a_j|_{p^1_n}  = x.
\end{equation}
\end{Lemma}

\begin{proof}
We apply Lemma \ref{technical} to the cohomology
classes $a_j$ and $\alpha_j$.
By  construction,
$a_j|_{p_0^1} =  0$ 
and $\sum_i a_j|_{p_1^i} = a_j|_{p_1^j} = x = \sum_i \alpha_j|_{q_1^i}$.
This gives us the necessary agreement on two coordinates.
\end{proof}

\begin{Lemma}\labell{aaLemma} 
For any $0 \leq j \leq n$ and $0 \leq k \leq n$,
$$\sum_i (a_j|_{p_k^i})^2 = {\frac{(n-1)!}{(k-1)!(n-k)!}} \ x^2, $$
where the sum is taken over all fixed points of index $2k$.
\end{Lemma}

\begin{proof}
We apply Lemma \ref{technical} to the cohomology
classes $(a_j)^2$ and $(\alpha_j)^2$.
By  construction,
$(a_j|_{p_0^1})^2 =  0 = (\alpha_j|_{q_0^1})^2 $ 
and $\sum_i (a_j|_{p_1^i})^2 = ( a_j|_{p_1^j})^2 = x^2 =
 \sum_i (\alpha_j|_{q_1^i})^2 $
Finally, by equation \ref{atop},
$(a_j|_{p^1_n})^2  = x^2$.
This gives us the necessary agreement on three coordinates.
\end{proof}

Combining  Lemmas \ref{aLemma} and \ref{aaLemma}, we see:

\begin{Lemma}\label{two}
The restriction $a_i|_F$ is equal to $0$ or $x$ for all $i$
and all fixed points $F$.
\end{Lemma}

In order to prove Proposition \ref{prop:4.5},
we will need to introduce some new cohomology
classes. 
r
Given a fixed point $F$ of index $2 k_F$, 
there exists a unique cohomology class $b_F$
such that $b_F|_F = x^{n-k_F}$, and
$b|_{F'} = 0$ for all other fixed points $F'$ whose index
is greater than or equal to $2k_F$.
(The proof of this function is essentially
identical to the proof of Theorem~\ref{getclasses},  
except in this case the Morse function is  $-\mu$.)

\begin{Lemma}
Let $F$ be any fixed point of index $2 k_F$,
and let $b_F \in H_{S^1}^{n-k_F}(M)$ be defined as above.
Then $a_j|_F \ x^{n-k_F} = b_F|_{p_1^j} \ x$.
\end{Lemma}

\begin{proof}

Pick any fixed point $\tilde{F}$ of index $2 k_F$ in $(P^1)^n$ such
that $\alpha_j|_{\tilde{F}} = a_j|_F$, and
let $\beta_{\tilde{F}} \in H_{S^1}^{n-k_F}(((P^1)^n),\Z)$ be the corresponding
cohomology class.
We  apply Lemma \ref{technical} to
$a_j \cdot b_F \in H^{2n-2k_F+2}(M,\Z)$ and $\alpha_j \cdot \beta_{\tilde{F}}$.
By construction,  $a_j$ vanishes on the unique fixed point
of index $0$, whereas $b_F$ vanishes on
all fixed points of index greater than $k_F$.
Thus $$\sum_i a_j|_{p^i_k} \cdot b_F|_{p^i_k} = 
\sum_i \alpha_j|_{q^i_k} \cdot \beta_{\tilde{F}}|_{q^i_k} = 0$$
for $k = 0$, or $k_F < k \leq n$.
Also, 
$$\sum_i a_j|_{p^i_{k_F}} \cdot b_F|_{p^i_{k_F}} = 
a_j|_F \cdot b_F|_F = a_j|_F = \alpha_j|_{\tilde{F}}$$
This gives the required agreement on $n-k_F+2$ coordinates.

Therefore
$$\sum_i a_j|_{p^i_1} \cdot b_F|_{p^i_1} = 
a_j|_{p^j_1} \cdot b_F|_{p^j_1} = 
\alpha_j|_{p^j_1} \cdot \beta_F|_{p^j_1},$$ 
and the result follows.

\end{proof}

Proposition \ref{prop:4.5} will  follow from the two lemmas below.

\begin{Lemma}\labell{lemma:4.13}
For each fixed point $F$ of index $2 k_F$, there
exist precisely $k_F$ numbers $j \in \{0,\ldots n\}$
such that $(a_j)_F = x$.
\end{Lemma}

\begin{proof}
Pick any fixed point $\tilde{F}$ in $(P^1)^n$ of index $2 k_F$. 
Let $b_F \in H_{S^1}^{2(n-k_F)}(M)$ 
and $\beta_{\tilde{F}} \in H_{S^1}^{2(n-k_F)}((\CP^1)^n)$ be the
classes  associated to $F$ and $\tilde{F}$, respectively.
We apply Lemma \ref{technical} to $b_F$  and $\beta_{\tilde{F}}$.
By construction,  both forms vanish on all points of index
greater than $2 k_F$, and 
the sum of each over points of index $2 k_F$ is $x^(n-k_F)$.
This gives agreement on $n - k_F + 1 $ coordinates.

Therefore
$\sum_i b_F|_{p_1^i}  =  k_F x^{n-k_F}$.
Combining this with the previous lemmas, the result follows.
\end{proof}

\begin{Lemma}
Let $J \subset \{1,\ldots,n\}$ be a subset with $|J|$ 
elements. There exists a unique fixed point $F$
of index $2{|J|}$ so that $a_j|_F = x$ for all $j \in J$.

\end{Lemma}

\begin{proof}
We apply Lemma \ref{technical} to
$a_J = \prod_{j \in J} a_j$ and $\alpha_J = \prod_{j \in J} \alpha_j$.
The restriction $\alpha_J|_F = 0$ for all critical points $F$ of index  
less that $2|J|$.  By the previous lemma, the same is true for $a_J$.
Also, $\alpha_J|_{q_n^1} = a_J|_{p_n^1} = x^{|J|}$.
This gives the necessary agreement on $|J| + 1$ coordinates.

Therefore, $\sum a_J|_{p_{|J|}^i} = \sum \alpha_J|_{q_{|J|}^i} = x^{|J|}$.
\end{proof}

\section{Cohomology rings of reduced spaces}

One corollary of Theorem 2 is that the cohomology of the reduced
spaces can also be computed.  Again the classical manifold invariants
are given by the same formulas as in the case of $(P^1)^n$.

\begin{Corollary}
Let $M$ be a compact, connected symplectic manifold equipped with
a semifree Hamiltonian circle action with isolated fixed points. 
Denote the moment map for the circle action by $\mu$, and 
suppose $0$ is a regular value of the moment map.
Then the cohomology ring of $M_{\rm red}$,
the reduced space at zero, is given by
$$H^*(M_{\rm red};\Z) = \Z[a_1,\ldots,a_n,y]/ R,$$
where $R$ is generated by the elements
\begin{enumerate}
\item $a_i y - a_i^2$ for all $i \in (1,\ldots, n)$, 
\item $ \prod_{j \in J} a_j$ for all $J \subset (1,\ldots,n)$ 
such that $\mu(F) >0$ for the corresponding fixed point $F$, and
\item  $ \prod_{j \not\in J} (y- a_j)$ for all $J \subset (1,\ldots,n)$
such that $\mu(F) < 0$ for the corresponding fixed point $F$.  
\end{enumerate}
The  Chern series of the reduced space
is given by
$$c_t = \prod_i 1+t(2a_i - y),$$
\end{Corollary}

For the case that $M = (P^1)^n$, the reduced space is a toric
variety, and the result above was proved by Haussmann and Knutson
\cite{HK}.
For the general case,  we use the following proposition from \cite{TW}:

\begin{Proposition} \labell{S1-Z}
Let $S^1$ act on a compact symplectic manifold $M$ with
moment map $\mu: M \to R$.  Assume that $0$ is a regular value of the 
moment map.
Let $F$ denote the set of fixed points.
Assume that
for every prime $p$, one of the following two conditions is satisfied:
\begin{enumerate}
\item  The integral cohomology of $F$ has no $p$-torsion, or:
\item For every point $m \in M$ which is not fixed by the $S^1-$action,
there exists a subgroup of $S^1$ congruent to $\Z/p$ which acts freely
on $m$.
\end{enumerate}

Define
$$K_+ := \{ \alpha \in H_{S^1}^*(M;\Z) \mid \alpha|_{F_+} = 0 \},
\mbox{    where }
F_+ := F \cap \phi\inv(0, \infty); \mbox{    }$$
$$K_- := \{ \alpha \in H_{S^1}^*(M;\Z) \mid \alpha|_{F_-} = 0 \},
\mbox{    where }
F_- := F \cap \phi\inv(-\infty, 0); \mbox{ and} $$
$$K := K_+ \oplus K_-.$$

Then there is a short exact sequence:
$$ 0 \to K \to H_{S^1}^*(M;\Z) \stackrel{\kappa}{\to}
H^*(M_\red;\Z) \to 0,$$
where $\kappa: H_{S^1}^*(M;\Z) \to H^*(\Mred;\Z)$ is  the Kirwan map.
\end{Proposition}

Note the case of a manifold equipped with a semi-free circle
action which has isolated fixed points is one where there is
no torsion at all in the cohomology of the fixed point set, and
where any $\Z/k$-subgroup of the circle acts freely outside fixed
points.  Hence both of the conditions in Proposition \ref{S1-Z}
are satisfied.

By another slight variation of Theorem \ref{getclasses},
for each fixed point $F$  of index
$2k$ such that $\mu(F) > 0$,
there exists a cohomology class $A_F \in K_+$ such that
$A_F|_F = x^{k}$ and so that $A_F|_{F'} = 0$ for all
other fixed points $F'$ such that  the index
of $F'$ is less than or equal to $k$,
Moreover, taken together over all fixed points,
there classes are a basis for $K_+$ as a module.
Because these classes are unique,
it is clear that $A_F =  \prod_{j \in J} a_j$,  
where $J$ is the corresponding subset of $\{1,\ldots,n\}$.

Similarly, $K_-$ is generated by
cohomology classes of the form
$ \prod_{j \not\in J} (y- a_j)$
for all $J \subset (1,\ldots,n)$ such that $\mu(F) < 0$ for the corresponding
fixed point $F$.

\end{document}